\documentclass[11pt]{article}
\usepackage{amsmath,amssymb,latexsym,epsf,mathabx}
\usepackage[all]{xy}
\usepackage{mathrsfs}
\usepackage{soul}
\usepackage{color}
\begin{document}

\newcommand{\nc}{\newcommand}
\newtheorem{Th}{Theorem}[section]
\newtheorem{Def}[Th]{Definition}
\newtheorem{Lem}[Th]{Lemma}
\newtheorem{Con}[Th]{Condition}
\newtheorem{Pro}[Th]{Proposition}
\newtheorem{Cor}[Th]{Corollary}
\newtheorem{Rem}[Th]{Remark}
\newtheorem{Exm}[Th]{Example}
\newtheorem{Sc}[Th]{}
\def\Pf#1{{\noindent\bf Proof}.\setcounter{equation}{0}}
\def\bskip#1{{ \vskip 20pt }\setcounter{equation}{0}}
\def\sskip#1{{ \vskip 5pt }\setcounter{equation}{0}}
\def\mskip#1{{ \vskip 10pt }\setcounter{equation}{0}}
\def\bg#1{\begin{#1}\setcounter{equation}{0}}
\def\ed#1{\end{#1}\setcounter{equation}{0}}

\soulregister\cite7 
\soulregister\citep7 
\soulregister\citet7 
\soulregister\ref7


\title{\bf  Ideal balanced pairs in $n$-exangulated categories
}
\smallskip
\author{ Peiyu Zhang, Zhuolin Zheng, Li Wang, Dajun Liu\\ 
}

\date{}
\maketitle
\baselineskip 15pt
%
%
\begin{abstract}
\vskip 10pt%

This paper is devoted to ideal approximation theory in $n$-exangulated categories. We define ideal balanced pairs and explore their properties and equivalent characterizations.
We prove that any complete ideal cotorsion triple $(\mathcal{I},\mathcal{K},\mathcal{J})$ satisfying suitable assumptions yields an ideal balanced pair
$(\mathcal{I}, \mathcal{J})$, and the converse holds under additional hypotheses. We further show that for an ideal balanced pair,
it is equivalent to an additive subfunctor with enough projective morphisms and enough injective morphisms.


\mskip\




\end{abstract}
%
\vskip 30pt

\section{Introduction}

Higher homological algebra has become a vibrant research field, whose core goal is to establish a higher-dimensional version of classical homological theories.
To construct a higher-dimensional Auslander correspondence and extend the classical theory of almost split sequences,
Iyama \cite{Iyama} put forward the concept of $n$-cluster tilting subcategories for arbitrary positive integers $n$.

Investigations on $n$-cluster tilting subcategories within abelian and exact categories motivated Jasso \cite{J} to propose new frameworks including $n$-abelian categories and $n$-exact categories.
In parallel work, Geiss, Keller and Oppermann \cite{GKO} defined $n$-angulated categories, which were developed to describe $(n-2)$-cluster tilting subcategories inside triangulated categories.
More recently, Herschend, Liu and Nakaoka \cite{HLN,HLN2} introduced the definition of $n$-exangulated categories for all positive integers $n$.
This class of categories serves as a higher generalization of the extriangulated categories introduced by Nakaoka and Palu \cite{NP}.
Meanwhile, it naturally unifies $n$-exact categories and $(n+2)$-angulated categories. Readers can refer to \cite{HHZL,HZhou,ZWei,ZZhou} for further research on these higher categorical structures.

The ideal approximation theory for an exact category \cite{FGHT}, as a generalization of approximation theory \cite{AR,GT},
aims to give a systematic study of the phenomena of approximations relative to an ideal.
It is devoted to the study of special precovering ideals and special preenveloping ideals \cite{FGHT,FH,FHH,O} and has produced many interesting applications.
For instance, Asadollahi and Sadeghi introduced ideal approximation theory into higher homological algebra \cite{AS},
and gave some interesting applications.
Wang and Wei \cite{WW} studied the ideal approximation theory associated with almost $n$-exact structures in the $n$-exangulated category.
The notions of $n$-ideal cotorsion pairs and $n$-$\mathbb{F}$-phantom morphisms were introduced and studied.
In particular, let $\mathscr{C}$ be an extriangulated category which satisfies the condition (WIC) and $\mathscr{T}$ be a nicely
embedded $n$-cluster tilting subcategory of $\mathscr{C}$, they proved Salce's Lemma in $\mathscr{T}$.

In a recent work, ideal approximation theory was introduced in extriangulated categories \cite{XFGS}.
The authors mainly introduced three lemmas in an extriangulated category. It is also shown that a finite intersection of special precovering (respectively, special preenveloping) ideals remains special precovering (respectively, special preenveloping). The results in this article show that extriangulated categories are the appropriate context for developing ideal approximation theory.

The concept of balanced pairs in abelian categories was originally introduced by Chen in \cite{CHEN}.
This concept is closely related to cotorsion pairs. The author proved that if the triple ($\mathcal{I}$, $\mathcal{K}$, $\mathcal{J}$) is a cotorsion triple
which is complete and hereditary, then the pair ($\mathcal{I}$, $\mathcal{J}$) is a balanced pair. In addition, he gave an equivalent characterization of balanced pairs.
A notable result concerning balanced pairs, obtained by Wang-Li-Huang \cite{WLH}, states that for an abelian category $\mathcal{A}$ whose $\mathrm{Ext}$ groups are small,
there exists a bijection between any two of the following three collections: (1) balanced pairs;
(2) subfunctors $\mathcal{F}\subseteq \mathrm{Ext}$ such that $\mathcal{A}$ admits enough $\mathcal{F}$-projectives and enough $\mathcal{F}$-injectives;
(3) Quillen exact structures $\mathcal{E}$ on $\mathcal{A}$ with enough $\mathcal{E}$-projectives and enough $\mathcal{E}$-injectives.

In this paper, we mainly study the ideal approximation theory, and investigate some properties and characterizations of ideal balanced pairs in the $n$-extriangulated category,
which is a higher-dimensional version of \cite{XF}. The rest of the paper is organized as follows.

In Section 2, we first recall some definitions and properties of the $n$-extriangulated category and ideals.
In Section 3, let $\mathcal{I}$ be a strongly precovering ideal and $C$ be an object in an $n$-extriangulated category. We introduce the notion of
an $\mathcal{I}$-resolution of $C$, and give the definition of an ideal balanced pair.
Given a complete ideal cotorsion triple $(\mathcal{I},\mathcal{K},\mathcal{J})$ satisfying certain conditions,
we show that it induces an ideal balanced pair $(\mathcal{I}, \mathcal{J})$ (see Theorem \ref{Th0}), with the converse holding under suitable hypotheses (see Theorem \ref{Th2}).
Meanwhile, we obtain two equivalent characterizations of an ideal balanced pair, (see Theorem \ref{Th1}).
In Section 4, we introduce the notion of an object ideal balanced pair, and establish the relationship between ideal balanced pairs $(\mathcal{I}, \mathcal{J})$
and $(\mathrm{Ob}\mathcal{I}, \mathrm{Ob}\mathcal{J})$. We prove that $(\mathcal{I},\mathcal{J})$ forms an ideal balanced pair
if and only if there exists an additive subfunctor $\mathbb{F}=\mathbb{E}_{\mathcal{I}}=\mathbb{E}^{\mathcal{J}}$ with enough projective morphisms and enough injective morphisms, and $\mathcal{I}=\mathbb{F}$-$\mathbf{proj}$, $\mathcal{J}=\mathbb{F}$-$\mathbf{inj}$.

\section{Preliminaries}
First, we recall some definitions and properties of $n$-extriangulated categories and ideals. Throughout this section, we assume that $\mathcal{C}$ is an additive category
and $n$ is a positive integer.
We denote by $\mathrm{Hom}_{\mathcal{C}}(A,B)$ or $\mathcal{C}(A,B)$ the set of morphisms from $A$ to $B$ in $\mathcal{C}$.
We use $\mathrm{Ab}$ to denote the category of abelian groups.

\subsection{$n$-Exangulated categories}
We recall some definitions and basic properties of $n$-exangulated categories. For more details, the reader can refer to \cite{HLN,HLN2}.
Suppose that $\mathcal{C}$ is equipped with an additive bifunctor $\mathbb{E}: \mathcal{C}^\mathrm{op}\times \mathcal{C}\to \mathrm{Ab}$. For any pair of objects $A,C\in \mathcal{C}$, an element $\delta\in \mathbb{E}(C,A)$ is called an $\mathbb{E}$-extension or simply an extension. We also write such $\delta$ as ${}_A\delta_C$ when we indicate $A$ and $C$. For any pair of $\mathbb{E}$-extensions ${}_A\delta_C$ and ${}_{A'}\delta'_{C'}$, let $\delta\oplus\delta'\in \mathbb{E}(C\oplus C',A\oplus A')$ be the element corresponding to $(\delta,0,0,\delta')$ through the natural isomorphism
\[
\mathbb{E}(C\oplus C',A\oplus A')\simeq \mathbb{E}(C,A)\oplus \mathbb{E}(C',A)\oplus \mathbb{E}(C,A')\oplus \mathbb{E}(C',A').
\]

For any $a\in \mathcal{C}(A,A')$ and $c\in \mathcal{C}(C,C')$, we have $\mathbb{E}(C,a)(\delta)\in \mathbb{E}(C,A')$ and $\mathbb{E}(c,A)(\delta)\in \mathbb{E}(C',A)$, which are denoted by $a_\ast\delta$ and $c^\ast\delta$, respectively.

Let ${}_A\delta_C$ and ${}_{A'}\delta'_{C'}$ be any pair of $\mathbb{E}$-extensions. A morphism of extensions $(a,c): \delta\to \delta'$ is a pair of morphisms $a\in \mathcal{C}(A,A')$ and $c\in \mathcal{C}(C,C')$ such that $a_\ast\delta=c^\ast\delta'$.

\medskip
\begin{Def}\label{}\cite[Definition 2.7]{HLN}.
Let $\mathrm{C}_\mathcal{C}$ be the category of complexes in $\mathcal{C}$. Define $\mathbf{C}_\mathcal{C}^{n+2}$ to be the full subcategory of $\mathbf{C}_\mathcal{C}$ consisting of the complexes whose components are zero in the degrees outside of $\{0,1,\dots,n+1\}$. Namely, an object in $\mathbf{C}_\mathcal{C}^{n+2}$ is a complex $X^\bullet=\{X,d_X^i\}$ of the form
$$\xymatrix{
X_{0}\ar[r]^{d_X^0}&X_{1}\ar[r]^{d_X^1}&...\ar[r]&X_{n}\ar[r]^{d_X^n}&X_{n+1}.
}$$
\end{Def}

We simply write a morphism $f^\bullet: X^\bullet\to Y^\bullet$ as $f^\bullet=(f^0,f^1,\dots,f^{n+1})$, only indicating the terms of degrees $0,\dots,n+1$.

By Yoneda's lemma, any $\mathbb{E}$-extension $\delta\in \mathbb{E}(C,A)$ induces natural transformations
\[
\delta_\sharp: \mathcal{C}(-,C)\Rightarrow \mathbb{E}(-,A)
\quad\text{and}\quad
\delta^\sharp: \mathcal{C}(A,-)\Rightarrow \mathbb{E}(C,-).
\]
For any $X\in \mathcal{C}$, these $(\delta_\sharp)_X$ and $(\delta^\sharp)_X$ are defined by
$(\delta_\sharp)_X: \mathcal{C}(X,C)\to \mathbb{E}(X,A),f\mapsto f^\ast\delta$
and
$(\delta^\sharp)_X: \mathcal{C}(A,X)\to \mathbb{E}(C,X),g\mapsto g_\ast\delta$.
In what follows, we may also simply write $\delta_X^\sharp(f)$ and $(\delta_\sharp)_X(g)$ as $\delta^\sharp(f)$ and $\delta_\sharp(g)$, respectively.

\medskip
\begin{Def}\label{}\cite[Definition 2.9]{HLN}.
Let $\mathcal{C},\mathbb{E},n$ be as before. Define the category $\AE:=\AE_\mathcal{C,\mathbb{E}}^{(n+2)}$ as follows.

$\mathrm{(1)}$ An object in $\AE_\mathcal{C,\mathbb{E}}^{(n+2)}$ is a pair $\langle X^\bullet,\delta\rangle$ consisting of $X^\bullet\in \mathbf{C}_\mathcal{C}^{n+2}$ and $\delta\in \mathbb{E}(X_{n+1},X_0)$ such that $(d_X^0)_\ast\delta=0$ and $(d_X^n)^\ast\delta=0.$
We call such a pair an $\mathbb{E}$-attached complex of length $n+2$. We also denote it by
$$\xymatrix{
X_{0}\ar[r]^{d_X^0}&X_{1}\ar[r]^{d_X^1}&...\ar[r]&X_{n}\ar[r]^{d_X^n}&X_{n+1}\ar@{.>}[r]^{\delta}&
}$$
$\mathrm{(2)}$ For any pairs $\langle X^\bullet,\delta\rangle$ and $\langle Y^\bullet,\rho\rangle$, a morphism $f^\bullet: \langle X^\bullet,\delta\rangle\to \langle Y^\bullet,\rho\rangle$ is defined to be a morphism $f^\bullet\in \mathbf{C}_\mathcal{C}^{n+2}(X^\bullet,Y^\bullet)$ satisfying $(f_0)_\ast\delta=(f_{n+1})^\ast\rho$.

We use the same composition and the identities as in $\mathbf{C}_\mathcal{C}^{n+2}$.
\end{Def}

\begin{Def}\label{}\cite[Definition 2.13]{HLN}.
An $n$-exangle is a pair $\langle X^\bullet,\delta\rangle$ consisting of $X^\bullet\in \mathbf{C}_\mathcal{C}^{n+2}$ and $\delta\in \mathbb{E}(X_{n+1},X_0)$ such that the following conditions hold.

$\mathrm{(1)}$ The following sequence of functors $\mathcal{C}^\mathrm{op}\to \mathrm{Ab}$ is exact:
$$\xymatrix{
\mathcal{C}(-,X_0)\ar[r]&\mathcal{C}(-,X_1)\ar[r]&\cdots\ar[r]&
\mathcal{C}(-,X_{n+1})\ar[r]^{\delta_\sharp}&\mathbb{E}(-,X_0)
}$$

$\mathrm{(2)}$ The following sequence of functors $\mathcal{C}\to \mathrm{Ab}$ is exact:
$$\xymatrix{
\mathcal{C}(X_{n+1},-)\ar[r]&\cdots\ar[r]&\mathcal{C}(X_0,-)\ar[r]^{\delta^\sharp}&\mathbb{E}(X_{n+1},-).
}$$
\end{Def}

In particular, any $n$-exangle is an object in $\mathbb{E}$. A morphism of $n$-exangles simply means a morphism in $\mathbb{E}$. Thus $n$-exangles form a full subcategory of $\mathbb{E}$.

\medskip
\begin{Pro}\label{claim2.15}\cite[Claim 2.15]{HLN}.
For any $n$-exangle $\langle X^\bullet,\delta\rangle$, the following are equivalent.

$\mathrm{(1)}$ $\delta=0$.

$\mathrm{(2)}$ There is $r\in \mathcal{C}(X_1,A)$ such that $r\circ d_X^0=1_A$.

$\mathrm{(3)}$ There is $s\in \mathcal{C}(C,X_n)$ such that $d_X^n\circ s=1_C$.

If one of the above conditions holds, then we say that the $n$-exangle $\langle X^\bullet,\delta\rangle$ splits.
\end{Pro}

\begin{Def}\label{}\cite[Definition 2.22]{HLN}.
Let $\mathfrak{s}$ be a correspondence which associates a homotopy equivalence class $\mathfrak{s}(\delta)=[{}_A X^\bullet_C]$ to each extension $\delta={}_A\delta_C$. Such $\mathfrak{s}$ is called a realization of $\mathbb{E}$ if it satisfies the following condition for any $\mathfrak{s}(\delta)=[X^\bullet]$ and $\mathfrak{s}(\rho)=[Y^\bullet]$.

$\mathrm{(R0)}$ For any morphism of extensions $(a,c): \delta\to \rho$, there exists a morphism $f^\bullet\in \mathbf{C}_\mathcal{C}^{n+2}(X^\bullet,Y^\bullet)$ of the form $f^\bullet=(a,f_1,\dots,f_n,c)$. Such $f^\bullet$ is called a lift of $(a,c)$.

In such a case, we simply say that ``$X^\bullet$ realizes $\delta$'' whenever they satisfy $\mathfrak{s}(\delta)=[X^\bullet]$. Moreover, a realization $\mathfrak{s}$ of $\mathbb{E}$ is said to be exact if it satisfies the following conditions.

$\mathrm{(R1)}$ For any $\mathfrak{s}(\delta)=[X^\bullet]$, the pair $\langle X^\bullet,\delta\rangle$ is an $n$-exangle.

$\mathrm{(R2)}$ For any $A\in \mathcal{C}$, the zero element ${}_A0_0=0\in \mathbb{E}(0,A)$ satisfies
$$\xymatrix{
\mathfrak{s}({}_A0_0)=[A\ar[r]^{1}&A\ar[r]&0\ar[r]&\cdots\ar[r]&0].
}$$
Dually, $\mathfrak{s}({}_00_A)=\bigl[0\to 0\to\cdots\to A\xrightarrow{1_A} A\bigr]$ holds for any $A\in \mathcal{C}$.
\end{Def}

\medskip
\begin{Def}\label{}\cite[Definition 2.23]{HLN}.
Let $\mathbb{E}$ and $\mathbf{C}_\mathcal{C}^{n+2}$ be as before, $\mathfrak{s}$ be an exact realization of $\mathbb{E}$. We use the following terminologies:

$\mathrm{(1)}$ An $n$-exangle $\langle X^\bullet,\delta\rangle$ is called an $\mathfrak{s}$-distinguished $n$-exangle if it satisfies $\mathfrak{s}(\delta)=[X^\bullet]$.

$\mathrm{(2)}$ An object $X^\bullet\in \mathbf{C}_\mathcal{C}^{n+2}$ is called an $\mathfrak{s}$-conflation or simply a conflation if it realizes some extension $\delta\in \mathbb{E}(X_{n+1},X_0)$.

$\mathrm{(3)}$ A morphism $f$ in $\mathcal{C}$ is called an $\mathfrak{s}$-inflation or simply an inflation if it admits some conflation $X^\bullet\in \mathbf{C}_\mathcal{C}^{n+2}$ satisfying $d_X^0=f$.

$\mathrm{(4)}$ A morphism $g$ in $\mathcal{C}$ is called an $\mathfrak{s}$-deflation or simply a deflation if it admits some conflation $X^\bullet\in \mathbf{C}_\mathcal{C}^{n+2}$ satisfying $d_X^n=g$.

When $\mathfrak{s}$ is clear from the context, we often simply say it is a distinguished $n$-exangle, conflation, inflation, or deflation.
\end{Def}

\begin{Def}\label{}\cite[Definition 2.27]{HLN}.
For a morphism $f^\bullet\in \mathbf{C}_\mathcal{C}^{n+2}(X^\bullet,Y^\bullet)$ satisfying $f_0=\mathrm{Id}_A$ for some $A=X_0=Y_0$, its mapping cone $M_f^\bullet\in \mathbf{C}_\mathcal{C}^{n+2}$ is defined to be the complex
$$\xymatrix{
X_{1}\ar[r]&X_{2}\bigoplus Y_{1}\ar[r]&X_{3}\bigoplus Y_{2}\ar[r]&...\ar[r]&X_{n+1}\bigoplus Y_{n}\ar[r]&Y_{n+1}.
}$$
The mapping cocone is defined dually, for morphisms $h^\bullet$ satisfying $h^{n+1}=\mathrm{Id}$ in $\mathrm{C}_\mathcal{C}^{n+2}$.
\end{Def}
Now we can define the $n$-exangulated category.

\begin{Def}\label{}\cite[Definition 2.32]{HLN}.
An $n$-exangulated category is a triplet $(\mathcal{C},\mathbb{E},\mathfrak{s})$ consisting of an additive category $\mathcal{C}$, a biadditive functor $\mathbb{E}: \mathcal{C}^\mathrm{op}\times \mathcal{C}\to \mathrm{Ab}$, and its exact realization $\mathfrak{s}$ such that the following conditions hold.

$\mathrm{(EA1)}$ Let $A\xrightarrow{f} B\xrightarrow{g} C$ be any sequence of morphisms in $\mathcal{C}$. If both $f$ and $g$ are inflations, then so is $g\circ f$. Dually, if $f$ and $g$ are deflations then so is $g\circ f$.

$\mathrm{(EA2)}$ For $\rho\in \mathbb{E}(D,A)$ and $c\in \mathcal{C}(C,D)$, let ${}_A\langle X^\bullet,c^\ast\rho\rangle_C$ and ${}_A\langle Y^\bullet,\rho\rangle_D$ be distinguished $n$-exangles. Then $(1_A,c)$ has a good lift $f^\bullet$, in the sense that its mapping cone gives a distinguished $n$-exangle $\langle M_f^\bullet,(d_X^0)_\ast\rho\rangle$.

$\mathrm{(EA2^\mathrm{op})}$ Dual of $(\mathrm{EA2})$.
\end{Def}
Unless otherwise specified, the category $(\mathcal{C},\mathbb{E},\mathfrak{s})$ in this paper always refers to an $n$-exangulated category.
The following lemma is frequently used.

\begin{Lem}\label{lemma3.12}\cite[Section~3]{HLN}.
For any distinguished $n$-exangle ${}_A\langle X^\bullet,\delta\rangle_C$, the following sequences are exact.
$$\xymatrix{
\mathcal{C}(-,A)\ar[r]&\cdots\ar[r]&
\mathcal{C}(-,C)\ar[r]&\mathbb{E}(-,A)\ar[r]&\mathbb{E}(-,X_{1})\ar[r]&\mathbb{E}(-,X_{2})
}$$
and
$$\xymatrix{
\mathcal{C}(C,-)\ar[r]&\cdots\ar[r]&\mathcal{C}(A,-)\ar[r]&\mathbb{E}(C,-)\ar[r]&\mathbb{E}(X_{n},-)\ar[r]&\mathbb{E}(X_{n-1},-).
}$$

\end{Lem}

\subsection{Ideals}
By an \textit{ideal} $\mathcal{I}$ in an $n$-exangulated category $\mathcal{C}$, we mean an additive subfunctor of $\mathrm{Hom}: \mathcal{C}^{\mathrm{op}} \times \mathcal{C} \to \mathbf{Ab}$. For an ordered pair $(A,B)$, there is a subgroup $\mathcal{I}(A,B)\subseteq \mathrm{Hom}(A,B)$ such that for any morphism $i\in \mathcal{I}(A,B)$, and any morphisms $f: X\to A$ and $g: B\to Y$ in $\mathcal{C}$, the composition $gif$ belongs to $\mathcal{I}(X,Y)$.

If $A\in \mathcal{C}$ is an object, we say that $A$ belongs to $\mathcal{I}$ if the identity morphism $1_A$ belongs to $\mathcal{I}(A,A)$, and we denote by $\mathrm{Ob}(\mathcal{I})\subseteq \mathcal{C}$ the full subcategory of objects of $\mathcal{I}$. This is an additive subcategory of $\mathcal{C}$.
In the other direction, for an additive subcategory $\mathcal{X}\subseteq \mathcal{C}$, there is an ideal $\mathcal{I}(\mathcal{X})$ generated by identity morphisms $1_X$, where $X\in \mathcal{X}$. Every morphism in this ideal factors through an object in $\mathcal{X}$. Such an ideal is called an \textit{object ideal}.

For a collection $\mathcal{M}$ of morphisms in $\mathcal{C}$, we define the notions of left and right orthogonal ideals by setting
$$\xymatrix{
\mathcal{M}^{\perp}:=\{g \mid \mathbb{E}(m, g)=0, \forall m \in \mathcal{M}\}
}$$
and
$$\xymatrix{
{}^{\perp}\mathcal{M}:=\{f \mid \mathbb{E}(f, m)=0, \forall m \in \mathcal{M}\},
}$$
respectively.

\begin{Con}[WIC]\cite[Condition 5.8]{NP}.
Let $(\mathcal{C}, \mathbb{E}, \mathfrak{s})$ be an n-extriangulated category. Consider the following conditions.

$\mathrm{(1)}$ Let $f \in \mathcal{C}(A, B)$, $g \in \mathcal{C}(B, C)$ be any composable pair of morphisms. If $g \circ f$ is an inflation, then so is $f$.

$\mathrm{(2)}$ Let $f \in \mathcal{C}(A, B)$, $g \in \mathcal{C}(B, C)$ be any composable pair of morphisms. If $g \circ f$ is a deflation, then so is $g$.

\end{Con}

\section{Ideal balanced pairs}

In this section, we introduce the notion of ideal balanced pair, and give some properties and equivalent characterizations.
First, we present several special classes of ideals, which are very useful in the sequel.

\begin{Def}\label{}
Assume that $\mathcal{I}$ is an ideal of $\mathcal{C}$.

$\mathrm{(1)}$ For any $C\in \mathcal{C}$, there is a morphism $i: I\to C$ in $\mathcal{I}$ such that for every morphism $i': I'\to C$ in $\mathcal{I}$, $i'$ factors through $i$,
i.e., we have the following commutative diagram,
$$\xymatrix{
&I'\ar[d]^{i'}\ar@{.>}[dl]_{\exists}\\
I\ar[r]^{i}&C.\\
}$$
Then the morphism i is called an $\mathcal{I}$-precover of C, and the ideal $\mathcal{I}$ is called a precovering ideal.

$\mathrm{(2)}$ For any $C\in \mathcal{C}$, if there exists the following commutative diagram:
$$\xymatrix{
A\ar[r]\ar[d]^{j}&X_{1}\ar[r]\ar[d]&...\ar[r]&X_{n}\ar[r]\ar[d]&C\ar@{.>}[r]^{\delta}\ar[d]&\\
A'\ar[r]&Y_{1}\ar[r]&...\ar[r]&Y_{n}\ar[r]^{i}&C\ar@{.>}[r]^{\eta}&
}$$
with $i\in \mathcal{I}$, $j\in \mathcal{I}^{\bot}$ and $\eta=j_{\ast}\delta$,
then the morphism i is called special $\mathcal{I}$-precover of C, and the ideal $\mathcal{I}$ is called a special precovering ideal.

$\mathrm{(3)}$ For any $C\in \mathcal{C}$, if there exist a distinguished n-exangle:
$$\xymatrix{
A\ar[r]&X_{1}\ar[r]&...\ar[r]&X_{n}\ar[r]^{i}&C\ar@{.>}[r]^{\delta}&
}$$
such that $i$ is an $\mathcal{I}$-precover of $C$, then the ideal $\mathcal{I}$ is called a strongly precovering ideal.

$\mathrm{(4)}$ For any $C\in \mathcal{C}$, if there exist a distinguished n-exangle:
$$\xymatrix{
A\ar[r]&X_{1}\ar[r]&...\ar[r]&X_{n}\ar[r]^{i}&C\ar@{.>}[r]^{\delta}&
}$$
such that $i\in\mathcal{I}$ and $A\in \mathcal{I}^{\bot}$, then the ideal $\mathcal{I}$ is called an object-special precovering ideal.

Dually, we can define $\mathcal{I}$-preenvelope, preenveloping ideal, special preenveloping ideal, strongly preenveloping ideal and object-special preenveloping ideal.
\end{Def}

From the definitions above, we readily obtain the following inclusion relations:
\begin{align*}
\text{object-special precovering ideal}
&\subseteq \text{special precovering ideal} \\
&\subseteq \text{strongly precovering ideal} \\
&\subseteq \text{precovering ideal.}
\end{align*}

Recall that a morphism $p$ is called projective if for any distinguished $n$-exangle $\delta$, there exists a morphism such that the following
diagram is commutative.
$$\xymatrix{
&&&&X\ar[d]^{p}\ar@{.>}[dl]&\\
A\ar[r]&X_{1}\ar[r]&...\ar[r]&X_{n}\ar[r]^{i}&C\ar@{.>}[r]^{\delta}&.
}$$
The category $\mathcal{C}$ has enough projective morphisms if for any $C\in \mathcal{C}$, there exists a distinguished $n$-exangle
$$\xymatrix{
A\ar[r]&X_{1}\ar[r]&...\ar[r]&X_{n}\ar[r]^{p}&C\ar@{.>}[r]^{\delta}&
}$$
such that $p$ is projective. Dually, we can define injective morphisms and $\mathcal{C}$ has enough injective morphisms.
We denote projective morphisms and injective morphisms in $\mathcal{C}$ by $\mathbb{E}$-$\mathbf{proj}$ and $\mathbb{E}$-$\mathbf{inj}$, respectively.

\mskip\
The following remark is mainly inspired by \cite[Remark~2.9]{ZZ1}.

\begin{Rem}\label{}
Assume that the category $\mathcal{C}$ has enough projective morphisms and enough injective morphisms, and satisfies the condition $\mathrm{(WIC)}$.
Let $\mathcal{I}$ be an ideal of $\mathcal{C}$, then

$\mathrm{(1)}$ $\mathcal{I}$ is a strongly precovering ideal if and only if $\mathcal{ I}$ is a precovering ideal and contains $\mathbb{E}$-$\mathbf{proj}$.

$\mathrm{(2)}$ $\mathcal{I}$ is a strongly preenveloping ideal if and only if $\mathcal{I}$ is a preenveloping ideal and contains $\mathbb{E}$-$\mathbf{inj}$.
\end{Rem}

\Pf. We only prove that (1) holds.

$(\Rightarrow)$ It is sufficient to prove that $\mathbb{E}$-$\mathbf{proj} \subseteq \mathcal{I}$.
For any projective morphism $p:X\to C$ in $\mathcal{C}$, since $\mathcal{I}$ is a strongly precovering ideal, we have the commutative diagram:
$$\xymatrix{
&&&&X\ar[d]^{p}\ar@{.>}[dl]_{f}&\\
A\ar[r]&X_{1}\ar[r]&...\ar[r]&X_{n}\ar[r]^{i}&C\ar@{.>}[r]^{\delta}&.
}$$
with $i$ $\mathcal{I}$-precover of $\mathcal{C}$. So $p=if$, and then $p\in \mathcal{I}$ since $\mathcal{I}$ is an ideal.
i.e., $\mathbb{E}$-$\mathbf{proj} \subseteq \mathcal{I}$.

$(\Leftarrow)$ For any object $C \in\mathcal{C}$,
since $\mathcal{ I}$ is a precovering ideal, there is an $\mathcal{I}$-precover $i: X\to C$.
We only prove that $i$ is a deflation. Since $\mathcal{C}$ has enough projective morphisms, there is a projective deflation $p:Y\to C$.
Thus we have the following commutative diagram:
$$\xymatrix{
&Y\ar[d]^{p}\ar@{.>}[dl]_{f}\\
X\ar[r]^{i}&C\\
}$$
as $\mathcal{ I}$ contains $\mathbb{E}$-$\mathbf{proj}$. So we have that $p=if$ is a deflation, since $\mathcal{C}$ satisfies the condition $\mathrm{(WIC)}$, $i$ is a deflation.
i.e., $\mathcal{I}$ is a strongly precovering ideal.
\ \hfill $\Box$

\begin{Def}\label{resolution}
Let $\mathcal{I}$ be a strongly precovering ideal of $\mathcal{C}$ and $C$ an object in $\mathcal{C}$.
An $\mathcal{I}$-resolution of $C$ is the following induced distinguished n-exangles:
$$\xymatrix{
K_{0}\ar[r]&Z_{01}\ar[r]&...\ar[r]&Z_{0,n-1}\ar[r]&X_{0}\ar[r]^{i_{0}}&C\ar@{.>}[r]^{\delta_{0}}&
}$$
$$\vdots$$
$$\xymatrix{
K_{m}\ar[r]&Z_{m1}\ar[r]&...\ar[r]&Z_{m,n-1}\ar[r]&X_{m}\ar[r]^{i_{m}}&K_{m-1}\ar@{.>}[r]^{\delta_{m}}&
}$$
$$\vdots$$
such that all $i_{\alpha}$ are an $\mathcal{I}$-precover for any $\alpha\geq0$.

Dually, if $\mathcal{J}$ is a strongly preenveloping ideal in $\mathcal{C}$, we can define the $\mathcal{J}$-coresolution of objects.
\end{Def}

With the above definitions, we now give the definition of ideal balanced pairs as follows.

\begin{Def}\label{BP}
Let $\mathcal{I}$ and $\mathcal{J}$ be two ideals of $\mathcal{C}$.
The ordered pair $(\mathcal{I},\mathcal{J})$ is called an ideal balanced pair in $\mathcal{C}$ if the following conditions hold:

$\mathrm{(BP1)}$ $\mathcal{I}$ is a strongly precovering ideal and $\mathcal{J}$ is a strongly preenveloping ideal.

$\mathrm{(BP2)}$ For each object $C\in\mathcal{C}$, there is a distinguished $n$-exangle $\delta_{n}$ induced by the $\mathcal{I}$-resolution of $C$ satisfying
$j_{\ast}\delta_{n}=0$ for any $j\in \mathcal{J}$.

$\mathrm{(BP3)}$ For each object $C\in\mathcal{C}$, there is a distinguished $n$-exangle $\theta_{n}$ induced by the $\mathcal{J}$-coresolution of $C$ satisfying
$i^{\ast}\theta_{n}=0$ for any $i\in \mathcal{I}$.
\end{Def}

\begin{Rem}\label{}
If the ordered pair $(\mathcal{I},\mathcal{J})$ is an ideal balanced pair in $\mathcal{C}$, then
$\mathbb{E}$-$\mathbf{proj}\subseteq\mathcal{I}$ and $\mathbb{E}$-$\mathbf{inj}\subseteq\mathcal{J}$.
\end{Rem}

\Pf. For any projective morphism $p$: $M\to N$, since $\mathcal{I}$ is a strongly precovering ideal, we have the following
commutative diagram:
$$\xymatrix{
&&&&M\ar[d]^{p}\ar[dl]_{j}&\\
A\ar[r]&X_{1}\ar[r]&...\ar[r]&X_{n}\ar[r]^{i}&N\ar@{.>}[r]^{\delta}&.
}$$
with $i$ being an $\mathcal{I}$-precover of $N$. Thus we know that $p=ij$, and then $p\in \mathcal{I}$ as $\mathcal{I}$ is an ideal.
i.e., $\mathbb{E}$-$\mathbf{proj}\subseteq\mathcal{I}$. Dually, we can prove that $\mathbb{E}$-$\mathbf{inj}\subseteq\mathcal{J}$.
\ \hfill $\Box$

\mskip\
The following lemma will be frequently used in the sequel.

\begin{Lem}\label{lemma1}
Let $\xymatrix{X_{0}\ar[r]^{d_{0}}&X_{1}\ar[r]^{d_{1}}&...\ar[r]&X_{n}\ar[r]^{d_{n}}&X_{n+1}\ar@{.>}[r]^{\delta}&}$ be a distinguished n-exangle.

$\mathrm{(1)}$ If $f: A\to X_{n+1}$ is in $\mathcal{C}$, then $f$ factors through $d_{n}$ if and only if $f^{\ast}\delta=0$.

$\mathrm{(2)}$ If $g: X_{0}\to B$ is in $\mathcal{C}$, then $g$ factors through $d_{0}$ if and only if $g_{\ast}\delta=0$.
\end{Lem}

\Pf, We only prove that (1) holds. The proof of (2) is similar.

By Lemma \ref{lemma3.12}, the following sequence is exact
$$\xymatrix{
\mathcal{C}(A,X_{n})\ar[r]^{\beta}&\mathcal{C}(A,X_{n+1})\ar[r]^{\alpha}&\mathbb{E}(A,X_{0}).
}$$
So we have that $\mathrm{Ker}\alpha=\mathrm{Im}\beta$ and $\alpha(f)=f^{\ast}\delta$.
Thus $f^{\ast}\delta=0$ if and only if $f\in \mathrm{Ker}\alpha$ if and only if $f\in \mathrm{Im}\beta$
if and only if $f$ factors through $d_{n}$.
\ \hfill $\Box$

\begin{Exm}\label{}
Assume that the category $\mathcal{C}$ has enough projective morphisms and enough injective morphisms.
Then $\mathbb{E}$-$\mathbf{proj}$ is a strongly precovering ideal and $\mathbb{E}$-$\mathbf{inj}$ is a strongly preenveloping ideal.
By Lemma \ref{lemma1}, the conditions $\mathrm{(BP2)}$ and $\mathrm{(BP3)}$ hold. Thus $(\mathbb{E}$-$\mathbf{proj}$, $\mathbb{E}$-$\mathbf{inj})$
is an ideal balanced pair.
\end{Exm}

Recall a subcategory $\mathcal{X}\subseteq \mathcal{C}$ is closed under extensions if for any distinguished $n$-exangle
$$\xymatrix{
X_{0}\ar[r]&X_{1}\ar[r]&...\ar[r]&X_{n}\ar[r]&X_{n+1}\ar@{.>}[r]^{\delta}&.
}$$
with $X_{0}$ and $X_{n+1}\in \mathcal{X}$, then $X_{i} \in \mathcal{X}$ for all $1\leq i\leq n$.
The subcategory $\mathcal{X}$ is closed under cones, if $X_{n}$ and $X_{n+1}$ are in $\mathcal{X}$, then $X_{0} \in \mathcal{X}$.
The subcategory $\mathcal{X}$ is closed under cocones, if $X_{0}$ and $X_{1}$ are in $\mathcal{X}$, then $X_{n+1} \in \mathcal{X}$.

An object $P$ is called projective if for any distinguished $n$-exangle $\delta$ and morphism $f$, we have the following commutative diagram.
$$\xymatrix{
&&&&X\ar[d]^{f}\ar@{.>}[dl]&\\
A\ar[r]&X_{1}\ar[r]&\cdots\ar[r]&X_{n}\ar[r]^{d_{n}}&C\ar@{.>}[r]^{\delta}&.
}$$
The category $\mathcal{C}$ has enough projective objects if for any $C\in \mathcal{C}$, there exists a distinguished $n$-exangle
$$\xymatrix{
C'\ar[r]&P_{1}\ar[r]&\cdots\ar[r]&P_{n}\ar[r]&C\ar@{.>}[r]^{\delta}&
}$$
such that all $P_{i}$ are projective. Dually, we can define injective objects and $\mathcal{C}$ has enough injective objects.
We denote projective objects and injective objects in $\mathcal{C}$ by $\mathbf{P}$ and $\mathbf{I}$, respectively.

\begin{Lem}\label{lemma2}
Let $\mathcal{I}$ be a special precovering ideal. If $\mathcal{I}^{\bot}$ is an object ideal, then $\mathcal{I}$ is an object-special precovering ideal.
\end{Lem}

\Pf. For any $C\in \mathcal{C}$, since $\mathcal{I}$ is a special precovering ideal, there exists the following commutative diagram:
$$\xymatrix{
A\ar[r]\ar[d]^{j}&X_{1}\ar[r]\ar[d]&...\ar[r]&X_{n}\ar[r]\ar[d]&C\ar@{.>}[r]^{\delta}\ar[d]&\\
A'\ar[r]&Y_{1}\ar[r]&...\ar[r]&Y_{n}\ar[r]^{i'}&C\ar@{.>}[r]^{\eta}&
}$$
with $i'\in \mathcal{I}$, $j\in \mathcal{I}^{\bot}$ and $\eta=j_{\ast}\delta$.
Since $\mathcal{I}^{\bot}$ is an object ideal, there are two morphisms $\alpha: A\to Y$ and $\beta:Y\to A'$ with $Y\in \mathcal{I}^{\bot}$ such that
$j=\beta\alpha$. Thus we have the following commutative diagram:
$$\xymatrix{
A\ar[r]\ar[d]^{\alpha}&X_{1}\ar[r]\ar[d]&...\ar[r]&X_{n}\ar[r]\ar[d]&C\ar@{.>}[r]^{\delta}\ar@{=}[d]&\\
Y\ar[r]\ar[d]^{\beta}&M_{1}\ar[r]\ar[d]&...\ar[r]&M_{n}\ar[r]^{i}\ar[d]^{\gamma}&C\ar@{.>}[r]^{\alpha_{\ast}\delta}\ar@{=}[d]&\\
A'\ar[r]&Y_{1}\ar[r]&...\ar[r]&Y_{n}\ar[r]^{i'}&C\ar@{.>}[r]^{\eta}&.
}$$
From the diagram above, we have that $i=i'\gamma$. i.e., $i\in \mathcal{I}$.
It follows from the second row in the diagram above that $\mathcal{I}$ an is object-special precovering ideal.
\ \hfill $\Box$

\mskip\
Recall that a pair $(\mathcal{I},\mathcal{J})$ of ideals of $\mathcal{C}$ is called an \textit{ ideal cotorsion pair}\cite[Definition 3.4]{WW} if
$\mathcal{I}={}^{\perp}\mathcal{J}$ and $\mathcal{J}=\mathcal{I}^{\perp}$. If $\mathcal{I}$ is a special precovering ideal and $\mathcal{J}$ is a special preenveloping ideal,
then the pair is called a complete cotorsion pair.

A subcategory $\mathcal{X}\subseteq \mathcal{C}$ is \textit{resolving} if it contains $\mathbf{P}$ and is closed under extensions and cones.
The notion of \textit{coresolving class} is defined dually. The following result shows that an ideal balanced pair can arise from an ideal cotorsion triple naturally,
which is also the first result in this section.

\begin{Th}\label{Th0}
Assume that $\mathcal{C}$ is an n-exangulated category with enough injective and projective objects.
Let $\mathcal{I},\mathcal{K},\mathcal{J}$ be ideals of $\mathcal{C}$ with $\mathcal{K}$ an object ideal and $\mathrm{Ob}(\mathcal{K})$ resolving-coresolving.
If both $(\mathcal{I},\mathcal{K})$ and $(\mathcal{K},\mathcal{J})$ are complete ideal cotorsion pairs in $\mathcal{C}$,
then $(\mathcal{I},\mathcal{J})$ is an ideal balanced pair in $\mathcal{C}$.
\end{Th}

\Pf. (BP1) Since $(\mathcal{I},\mathcal{K})$ and $(\mathcal{K},\mathcal{J})$ are complete ideal cotorsion pairs,
$\mathcal{I}$ is a special precovering ideal and is also a strongly precovering ideal.
Similarly, $\mathcal{J}$ is a strongly preenveloping ideal.

(BP2) Since $(\mathcal{I},\mathcal{K})$ is a complete ideal cotorsion pair and $\mathcal{K}$ is an object ideal, we have
that $\mathcal{I}$ is an object-special precovering ideal and is also a strongly precovering ideal by Lemma \ref{lemma2}.
For any object $C\in \mathcal{C}$, each induced distinguished $n$-exangle (see Definition \ref{resolution}):
$$\xymatrix{
K_{m}\ar[r]&Z_{m1}\ar[r]&...\ar[r]&Z_{m,n-1}\ar[r]&X_{m}\ar[r]^{i_{m}}&K_{m-1}\ar@{.>}[r]^{\delta_{m}}&
}$$
satisfying $i_{m}$ is an $\mathcal{I}$-precover and $K_{m} \in \mathcal{I}^{\bot}=\mathcal{K}$.
It is sufficient to prove that $j_{\ast}\delta_{m}=0$ for any $j\in\mathcal{J}$.
Since $\mathcal{C}$ has enough injective objects, by Proposition 3.6 in \cite{HLN}, we have the following commutative diagram:
$$\xymatrix{
K_{m}\ar[r]\ar@{=}[d]&Z_{m1}\ar[r]\ar[d]&...\ar[r]&Z_{m,n-1}\ar[r]\ar[d]&X_{m}\ar[r]\ar[d]&K_{m-1}\ar@{.>}[r]^{\delta_{m}}\ar[d]^{f}&\\
K_{m}\ar[r]&E_{1}\ar[r]&...\ar[r]&E_{n-1}\ar[r]&E_{n}\ar[r]&K\ar@{.>}[r]^{\delta}&
}$$
with $E_{i}\in \mathbf{I}$ and $\delta_{m}=f^{\ast}\delta$. Since $\mathrm{Ob}(\mathcal{K})$ is coresolving, we have that $K\in \mathcal{K}$.
Note that $\mathbb{E}(K,j)(\delta)=j^{\ast}\delta=0$ as $(\mathcal{K},\mathcal{J})$ is an ideal cotorsion pair.
Thus $j_{\ast}\delta_{m}=j_{\ast}f^{\ast}\delta=f^{\ast}j_{\ast}\delta=0$.

(BP3) It is duall to (BP2).
\ \hfill $\Box$

\mskip\
Let $\mathcal{I}$ and $\mathcal{J}$ be two ideals of $\mathcal{C}$, and $A$ and $C$ be two objects in $\mathcal{C}$. We define two subclasses of $\mathbb{E}(C,A)$ as follows:
$$\mathbb{E}_{\mathcal{I}}(C,A)=\{\delta\in\mathbb{E}(C,A)\mid i^{*}\delta=0 \text{ for any }i\in\mathcal{I}\},$$
$$\mathbb{E}^{\mathcal{J}}(C,A)=\{\delta\in\mathbb{E}(C,A)\mid j_{*}\delta=0 \text{ for any }j\in\mathcal{J}\}.$$

Let $\mathbb{E}_{\mathcal{I}}$ be the class consisting of all distinguished $n$-exangles in $\mathbb{E}_{\mathcal{I}}(C,A)$, and
let $\mathbb{E}^{\mathcal{J}}$ be the class consisting of all distinguished $n$-exangles in $\mathbb{E}^{\mathcal{J}}(C,A)$ respectively,
where $A$ and $C$ are any two objects in $\mathcal{C}$.

\begin{Pro}\label{BP1}
Let $\mathcal{I}$ be a strongly precovering ideal, and $\mathcal{J}$ a strongly preenveloping ideal of $\mathcal{C}$.
Then the pair $(\mathcal{I},\mathcal{J})$ is an ideal balanced pair in $\mathcal{C}$ if and only if $\mathbb{E}_{\mathcal{I}}=\mathbb{E}^{\mathcal{J}}$.
\end{Pro}

\Pf. $(\Leftarrow)$ (BP1) is obvious and (BP3) is duall to (BP2). Next, we only need to prove that (BP2) holds.

For any object $C\in \mathcal{C}$, each induced distinguished $n$-exangle (see Definition \ref{resolution}):
$$\xymatrix{
K_{m}\ar[r]&Z_{m1}\ar[r]&...\ar[r]&Z_{m,n-1}\ar[r]&X_{m}\ar[r]^{i_{m}}&K_{m-1}\ar@{.>}[r]^{\delta_{m}}&
}$$
satisfying $i_{m}$ is an $\mathcal{I}$-precover. For any $i\in \mathcal{I}$, $i$ factors through $i_{m}$, by Lemma \ref{lemma1},
$i^{\ast}\delta_{m}=0$. Since $\mathbb{E}_{\mathcal{I}}=\mathbb{E}^{\mathcal{J}}$, $j_{\ast}\delta_{m}=i^{\ast}\delta_{m}=0$.
So the condition (BP2) holds.

$(\Rightarrow)$ For any distinguished $n$-exangle
$$\xymatrix{
A\ar[r]&Y_{1}\ar[r]&...\ar[r]&Y_{n}\ar[r]^{d_{n}}&C\ar@{.>}[r]^{\delta}&
}$$
in $\mathbb{E}_{\mathcal{I}}(C,A)$. Since $\mathcal{I}$ is a strongly precovering ideal, we consider the first induced distinguished $n$-exangle
in the $\mathcal{I}$-resolution of $C$:
$$\xymatrix{
K_{0}\ar[r]&Z_{01}\ar[r]&...\ar[r]&Z_{0,n-1}\ar[r]&X_{0}\ar[r]^{i_{0}}&C\ar@{.>}[r]^{\delta_{0}}&
}$$
with $i_{0}$ being $\mathcal{I}$-precover. Since $i_{0}^{\ast}\delta=0$, by Lemma \ref{lemma1} and the dual of Proposition 3.6 in \cite{HLN},
there is a morphism $b$ such that the following diagram is commutative:
$$\xymatrix{
K_{0}\ar[r]\ar[d]^{a}&Z_{01}\ar[d]\ar[r]&...\ar[r]&Z_{0,n-1}\ar[d]\ar[r]&X_{0}\ar[d]^{b}\ar[r]^{i_{0}}&C\ar@{=}[d]\ar@{.>}[r]^{\delta_{0}}&\\
A\ar[r]&Y_{1}\ar[r]&...\ar[r]&Y_{n-1}\ar[r]&Y_{n}\ar[r]^{d_{n}}&C\ar@{.>}[r]^{\delta}&.
}$$
Then $\delta=a_{\ast}\delta_{0}$. For any $j\in \mathcal{J}$, $j_{\ast}\delta=j_{\ast}a_{\ast}\delta_{0}=(ja)_{\ast}\delta_{0}$.
Note that $ja\in \mathcal{J}$ as $\mathcal{J}$ is a ideal. By the definition of balanced pairs,
we have that $j_{\ast}\delta=(ja)_{\ast}\delta_{0}=0$. i.e., $\delta\in\mathbb{E}^{\mathcal{J}}$, and then $\mathbb{E}_{\mathcal{I}}\subseteq\mathbb{E}^{\mathcal{J}}$.
Similarly, we can prove that $\mathbb{E}^{\mathcal{J}}\subseteq\mathbb{E}_{\mathcal{I}}$.
\ \hfill $\Box$

\mskip\
We now present equivalent characterizations of an ideal balanced pair.

\begin{Th}\label{Th1}
Let $\mathcal{I}$ be a strongly precovering ideal, and $\mathcal{J}$ a strongly preenveloping ideal of $\mathcal{C}$.
Then the following statements are equivalent:

$\mathrm{(1)}$ the pair $(\mathcal{I},\mathcal{J})$ is an ideal balanced pair in $\mathcal{C}$.

$\mathrm{(2)}$ $\mathbb{E}_{\mathcal{I}}=\mathbb{E}^{\mathcal{J}}$.

$\mathrm{(3)}$ For each object $C\in \mathcal{C}$, there are two distinguished $n$-exangles
$$\xymatrix{
A\ar[r]&X_{1}\ar[r]&...\ar[r]&X_{n}\ar[r]^{i}&C\ar@{.>}[r]^{\delta_{1}}&
}$$
and
$$\xymatrix{
C\ar[r]^{j}&Y_{1}\ar[r]&...\ar[r]&Y_{n}\ar[r]&B\ar@{.>}[r]^{\delta_{2}}&
}$$
with $i\in \mathcal{I}$, $j\in \mathcal{J}$ and $\delta_{1}$, $\delta_{2}\in\mathbb{E}_{\mathcal{I}}\bigcap\mathbb{E}^{\mathcal{J}}$.
\end{Th}

\Pf. (1)$\Rightarrow$(2) By Proposition \ref{BP1}.

(2)$\Rightarrow$(3) Since $\mathcal{I}$ is a strongly precovering ideal, for any $C\in \mathcal{C}$,
there is a distinguished $n$-exangle:
$$\xymatrix{
A\ar[r]&X_{1}\ar[r]&...\ar[r]&X_{n}\ar[r]^{i}&C\ar@{.>}[r]^{\delta_{1}}&
}$$
with $i$ being an $\mathcal{I}$-precover of $C$. By Lemma \ref{lemma1}, we have that $i'^{\ast}\delta_{1}=0$ for any $i'\in \mathcal{I}$.
i.e., $\delta_{1}\in\mathbb{E}_{\mathcal{I}}=\mathbb{E}^{\mathcal{J}}$.
Thus we obtain the first distinguished $n$-exangle. Dually, we can get the second distinguished $n$-exangle.

(3)$\Rightarrow$(1) For any $C\in \mathcal{C}$, by (3), we have the following distinguished $n$-exangle
$$\xymatrix{
K_{0}\ar[r]&Z_{01}\ar[r]&...\ar[r]&Z_{0,n-1}\ar[r]&X_{0}\ar[r]^{i_{0}}&C\ar@{.>}[r]^{\delta_{0}}&
}$$
with $i_{0}$ being an $\mathcal{I}$-precover and $j_{\ast}\delta_{0}=0$ for any $j\in \mathcal{J}$.
Continuing this process, we obtain an $\mathcal{I}$-resolution of $C$ satisfying the condition (BP2).
Dually, we can get a $\mathcal{J}$-coresolution of $C$ satisfying the condition (BP3) by the second distinguished $n$-exangle in (3).
\ \hfill $\Box$

\mskip\
The following conclusion indicates that one subcategory in an ideal balanced pair is completely determined by the other.

\begin{Cor}\label{}
Let $\mathcal{C}$ be an n-exangulated category. Then the following statements hold:

$\mathrm{(1)}$ If $(\mathcal{I}_{1},\mathcal{J})$ and $(\mathcal{I}_{2},\mathcal{J})$ are two ideal balanced pairs in $\mathcal{C}$,
then $\mathcal{I}_{1}=\mathcal{I}_{2}$.

$\mathrm{(2)}$ If $(\mathcal{I},\mathcal{J}_{1})$ and $(\mathcal{I},\mathcal{J}_{2})$ are two ideal balanced pairs in $\mathcal{C}$,
then $\mathcal{J}_{1}=\mathcal{J}_{2}$.
\end{Cor}

\Pf. Statements (1) and (2) are dual to each other, hence, we only prove that (1) holds.

For any $f$: $M\to C$ in $\mathcal{I}_{1}$, since $\mathcal{I}_{2}$ is a strongly precovering ideal,
there is a distinguished $n$-exangle:
$$\xymatrix{
A\ar[r]&X_{1}\ar[r]&...\ar[r]&X_{n}\ar[r]^{i_{2}}&C\ar@{.>}[r]^{\delta}&
}$$
with $i_{2}$ being an $\mathcal{I}_{2}$-precover of $C$. By Lemma \ref{lemma1}, we have that $x^{\ast}\delta=0$ for any $x\in \mathcal{I}_{2}$.
Thus $\delta\in\mathbb{E}_{\mathcal{I}_{2}}=\mathbb{E}^{\mathcal{J}}=\mathbb{E}_{\mathcal{I}_{1}}$, and then $f^{\ast}\delta=0$.
By Lemma \ref{lemma1} again, there is a morphism $g: M\to X_{n}$ such that $f=i_{2}g$.
So we have that $f\in \mathcal{I}_{2}$ as $\mathcal{I}_{2}$ is an ideal. i.e., $\mathcal{I}_{1}\subseteq\mathcal{I}_{2}$.
Similarly, we can prove that $\mathcal{I}_{2}\subseteq\mathcal{I}_{1}$.
\ \hfill $\Box$

\begin{Def}\label{}
Assume that $\mathcal{C}$ is an $n$-exangulated category and $\mathcal{I}$ is an ideal.
Let the following diagram be a morphism of distinguished $n$-exangles.
$$\xymatrix{
X_{0}\ar[r]\ar[d]^{f_{0}}&X_{1}\ar[r]\ar[d]^{f_{1}}&...\ar[r]&X_{n}\ar[r]\ar[d]^{f_{n}}&X_{n+1}\ar[d]^{f_{n+1}}\ar@{.>}[r]&\\
Y_{0}\ar[r]&Y_{1}\ar[r]&...\ar[r]&Y_{n}\ar[r]&Y_{n+1}\ar@{.>}[r]&.
}$$
We say that $\mathcal{I}$ is closed under extensions if $f_{0}$, $f_{n+1}\in\mathcal{I}$ implies $f_{i}\in\mathcal{I}$ for all $1\leq i\leq n$.
\end{Def}

The following result is the converse of Theorem \ref{Th0}.

\begin{Th}\label{Th2}
Assume that $\mathcal{C}$ is an $n$-exangulated category with enough projective objects and enough injective objects.
Let $\mathcal{I}$ and $\mathcal{J}$ be two object ideals satisfying the following properties:

$\mathrm{(1)}$ $\mathcal{I}$ is object-special precovering and $\mathcal{J}$ is object-special preenveloping.

$\mathrm{(2)}$ $\mathcal{I}^{\perp}$ and ${}^{\perp}\mathcal{J}$ are closed under extensions.

$\mathrm{(3)}$ $\mathcal{I}^{\perp}\bigcap\mathcal{I}\subseteq{}^{\perp}\mathcal{J}$ and ${}^{\perp}\mathcal{J}\bigcap\mathcal{J}\subseteq\mathcal{I}^{\perp}$.

$\mathrm{(4)}$ $(\mathcal{I},\mathcal{J})$ is an ideal balanced pair.

Then, there are two complete ideal cotorsion pairs $(\mathcal{I},\mathcal{K})$ and $(\mathcal{K},\mathcal{J})$.
\end{Th}

\Pf. By Theorem 3.8 and its dual in \cite{WW}, we have that $(\mathcal{I},\mathcal{I}^{\bot})$ and $(^{\bot}\mathcal{J},\mathcal{J})$ are two ideal
cotorsion pairs. Hence it suffices to show that $\mathcal{I}^{\bot}=$ $^{\bot}\mathcal{J}$. Setting $\mathcal{K}=\mathcal{I}^{\bot}=$ $^{\bot}\mathcal{J}$,
we complete the proof.

For any $f$: $A\to B$ in $\mathcal{I}^{\bot}$, since $\mathcal{I}$ is object-special precovering, there is a distinguished $n$-exangle:
$$\xymatrix{
K\ar[r]&X_{1}\ar[r]&...\ar[r]&X_{n}\ar[r]^{i}&B\ar@{.>}[r]^{\delta}&
}$$
with $i$ being an $\mathcal{I}$-precover of $B$ and $K\in \mathcal{I}^{\bot}$.
We consider the following commutative diagram:
$$\xymatrix{
K\ar[r]^{t_{0}}\ar@{=}[d]&Y_{1}\ar[r]\ar[d]^{h_{1}}&...\ar[r]&Y_{n}\ar[r]^{t_{n}}\ar[d]^{h_{n}}&A\ar[d]^{f}\ar@{.>}[r]^{f^{\ast}\delta}&\\
K\ar[r]&X_{1}\ar[r]&...\ar[r]&X_{n}\ar[r]^{i}&B\ar@{.>}[r]^{\delta}&.
}$$
By Lemma \ref{lemma1}, for any $x\in \mathcal{I}$, $x^{\ast}\delta=0$. i.e., $\delta\in \mathbb{E}_{\mathcal{I}}$.
By Theorem \ref{Th1}, $\delta\in \mathbb{E}^{\mathcal{J}}$. For any $j\in \mathcal{J}$, $j_{\ast}f^{\ast}\delta=f^{\ast}j_{\ast}\delta=0$.
i.e., $f^{\ast}\delta\in\mathbb{E}^{\mathcal{J}}$. Note that $1_{K}$ and $f$ are in $\mathcal{I}^{\bot}$, then $h_{1}$, ..., $h_{n}$ are in $\mathcal{I}^{\bot}$
since $\mathcal{I}^{\perp}$ is closed under extensions. As $\mathcal{I}$ and $\mathcal{I}^{\bot}$ are ideals, we have that
$ih_{n}\in \mathcal{I}^{\perp}\bigcap\mathcal{I}\subseteq{}^{\perp}\mathcal{J}$.

Let $g$: $M\to J$ be a morphism with $J\in \mathrm{Ob}(\mathcal{J})$, and then $g\in \mathcal{J}$ since $\mathcal{J}$ is an object ideal.
Next, we will show that $\mathbb{E}(f,g)=0$. We consider the following commutative diagram:
$$\xymatrix{
& \mathbb{E}(B,M) \ar[rr] \ar@{.>}[dd]^{\mathbb{E}(f,M)} \ar[ld]_{\mathbb{E}(B,g)}
&& \mathbb{E}(X_{n},M) \ar[ld]\ar[rr] \ar@{.>}[dd]&& \mathbb{E}(X_{n-1},M) \ar[dd]\ar[ld] \\
\mathbb{E}(B,J) \ar[rr] \ar[dd]_{\mathbb{E}(f,J)}&& \mathbb{E}(X_{n},J) \ar[rr] \ar[dd]&& \mathbb{E}(X_{n-1},J) \ar[dd]\\
& \mathbb{E}(A,M) \ar@{.>}[rr]^{\mathbb{E}(t_{n},M)} \ar@{.>}[ld]_{\mathbb{E}(A,g)}&& \mathbb{E}(Y_{n},M) \ar@{.>}[rr] \ar@{.>}[ld]_{\mathbb{E}(Y_{n},g)}&& \mathbb{E}(Y_{n-1},M)\ar[ld], \\
\mathbb{E}(A,J) \ar[rr]_{\mathbb{E}(t_{n},J)}&& \mathbb{E}(Y_{n},J) \ar[rr]&& \mathbb{E}(Y_{n-1},J).
}$$
Consider the left-hand cube, we obtain the following equations:
\begin{align*}
\mathbb{E}(t_{n},J)\mathbb{E}(A,g)\mathbb{E}(f,M)&=\mathbb{E}(Y_{n},g)\mathbb{E}(t_{n},M)\mathbb{E}(f,M)\\
&=\mathbb{E}(Y_{n},g)\mathbb{E}(ft_{n},M)\\
&=\mathbb{E}(Y_{n},g)\mathbb{E}(ih_{n},M)\\
&=\mathbb{E}(ih_{n},g)\\
&=0.
\end{align*}
By Lemma \ref{lemma3.12}, we consider the following exact sequence
$$\xymatrix{
\mathcal{C}(Y_{1},J)\ar[r]^{\mathcal{C}(t_{0},J)}&\mathcal{C}(K,J)\ar[r]&\mathbb{E}(A,J)\ar[r]^{\mathbb{E}(t_{n},J)}&\mathbb{E}(Y_{n},J).
}$$
For any $h$: $K\to J$, we have the following commutative diagram:
$$\xymatrix{
K\ar[r]^{t_{0}}\ar[d]_{h}&Y_{1}\ar[r]\ar[d]^{h_{1}}&...\ar[r]&Y_{n}\ar[r]^{t_{n}}\ar[d]^{h_{n}}&A\ar@{=}[d]\ar@{.>}[r]^{f^{\ast}\delta}&\\
J\ar[r]^{x}&Z_{1}\ar[r]&...\ar[r]&Z_{n}\ar[r]^{i}&A\ar@{.>}[r]^{h_{\ast}f^{\ast}\delta}&.
}$$
Note that $f^{\ast}\delta\in\mathbb{E}^{\mathcal{J}}$. Then $h_{\ast}f^{\ast}\delta=0$ since  $h$ is in $\mathcal{J}$.
By Proposition \ref{claim2.15}, $x$ is a section.
i.e., $h$ factors through $t_{0}$. Thus $\mathcal{C}(t_{0},J)$ is epic, and then $\mathbb{E}(t_{n},J)$ is mono.
So we have that $\mathbb{E}(f,g)=\mathbb{E}(A,g)\mathbb{E}(f,M)=0$.

For any $y$: $M\to N$ in $\mathcal{J}$, since $\mathcal{J}$ is an object ideal, there are two morphisms $y_{1}$: $M\to L$ and $y_{2}$: $L\to N$ satisfying
$L\in \mathcal{J}$ and $y=y_{2}y_{1}$. By the argument above, $\mathbb{E}(f,y_{1})=0$.
Thus $\mathbb{E}(f,y)=\mathbb{E}(f,y_{2}y_{1})=\mathbb{E}(f,y_{2})\mathbb{E}(f,y_{1})=0.$
So we can get that $f\in ^{\perp}\mathcal{J}$. i.e., $\mathcal{I}^{\bot}\subseteq {}^{\bot}\mathcal{J}$.
Dually, we can obtain that ${}^{\bot}\mathcal{J}\subseteq\mathcal{I}^{\bot}$. Consequently, $\mathcal{I}^{\bot}= {}^{\bot}\mathcal{J}$.
\ \hfill $\Box$

\section{Object ideal balanced pairs and additive subfunctors}

In this section, we introduce the notion of an object ideal balanced pair, and establish the relationship between ideal balanced pairs $(\mathcal{I}, \mathcal{J})$
and $(\mathrm{Ob}(\mathcal{I}), \mathrm{Ob}(\mathcal{J}))$. We mainly show that $(\mathcal{I},\mathcal{J})$ forms an ideal balanced pair
if and only if there exists an additive subfunctor satisfying some conditions.

\subsection{Object ideal balanced pairs}

\begin{Def}\label{}
An ideal balanced pair $(\mathcal{I},\mathcal{J})$ in $\mathcal{C}$ is called an object ideal balanced pair if both $\mathcal{I}$ and $\mathcal{J}$ are object ideals.
\end{Def}

\begin{Pro}\label{}
Let $\mathcal{I}$ and $\mathcal{J}$ be two object ideals. Then the following results hold.

$\mathrm{(1)}$ If $\mathcal{I}$ is a strongly precovering ideal, then $\mathrm{Ob}(\mathcal{I})$ is contravariantly finite in $\mathcal{C}$.

$\mathrm{(2)}$ If $\mathcal{J}$ is a strongly preenveloping ideal, then $\mathrm{Ob}(\mathcal{J})$ is covariantly finite in $\mathcal{C}$.
\end{Pro}

\Pf. For any object $C\in \mathcal{C}$, since $\mathcal{I}$ is an strongly precovering ideal, there is a deflation $i$: $X\to C$ such that $i$ is an $\mathcal{I}$-precover.
As $\mathcal{I}$ is an object ideal, there are two morphisms $i_{1}$: $X\to I$ and $i_{2}$: $I\to C$ satisfying $I\in \mathrm{Ob}(\mathcal{I})$ and $i=i_{2}i_{1}$.
Note that $i$ is an $\mathcal{I}$-precover. It is easy to show that the morphism $i_{2}$ is an $\mathrm{Ob}(\mathcal{I})$-precover of $C$.
i.e., we have the following commutative diagram:
$$\xymatrix{
&&\forall I'\in \mathrm{Ob}(\mathcal{I})\ar[d]^{\forall f}\ar@{.>}[dll]_{\exists g}\\
X\ar[rr]^{i}\ar[dr]_{i_{1}}&&C.\\
&I\ar[ur]_{i_{2}}&
}$$
Dually, we can show that $\mathrm{Ob}(\mathcal{J})$ is a covariantly finite subcategory in $\mathcal{C}$.
\ \hfill $\Box$

\begin{Def}\label{}
Let $\mathcal{C}$ be an $n$-extriangulated category. A subcategory $\mathcal{X}$ of $\mathcal{C}$ is called strongly exact-contravariantly finite,
if for any object $C\in\mathcal{C}$, there exists a distinguished $n$-exangle:
$$\xymatrix{
A\ar[r]&X_{1}\ar[r]&...\ar[r]&X_{n}\ar[r]&C\ar@{.>}[r]^{\delta}&
}$$
such that $X_{n} \in \mathcal{X}$ and the following sequence is exact:
$$\xymatrix{
0\ar[r]&\mathcal{C}(X,A)\ar[r]&\mathcal{C}(X,X_{1})\ar[r]&...\ar[r]&\mathcal{C}(X,X_{n})\ar[r]&\mathcal{C}(X,C)\ar[r]&0
}$$
for any $X\in \mathcal{X}$. i.e., the above distinguished $n$-exangle is $\mathcal{C}(\mathcal{X},-)$-exact.
The notion of a strongly exact-covariantly finite subcategory is defined dually.
\end{Def}

If $\mathcal{X}$ is strongly exact-contravariantly finite and $\mathcal{Y}$ is strongly exact-covariantly finite in $\mathcal{C}$,
similar to Definition \ref{resolution}, we can define the notions of $\mathcal{X}$-resolutions and $\mathcal{Y}$-coresolutions.

\begin{Def}\label{OBBP}
Let $\mathcal{X}$ and $\mathcal{Y}$ be two subcategories of $\mathcal{C}$ which are closed under direct summands.
The pair $(\mathcal{X},\mathcal{Y})$ is called a balanced pair in $\mathcal{C}$ if the following conditions hold:

$(1)$ $\mathcal{X}$ is strongly exact-contravariantly finite and $\mathcal{Y}$ is strongly exact-covariantly finite in $\mathcal{C}$.

$(2)$ For each object $C\in\mathcal{C}$, there is a distinguished $n$-exangle $\delta_{n}$ induced by the $\mathcal{X}$-resolution of $C$
such that $\delta_{n}$ is $\mathcal{C}(-,\mathcal{Y})$-exact.

$(3)$ For each object $C\in\mathcal{C}$, there is a distinguished $n$-exangle $\theta_{n}$ induced by the $\mathcal{Y}$-coresolution of $C$
such that $\theta_{n}$ is $\mathcal{C}(\mathcal{X},-)$-exact..
\end{Def}

The following proposition follows from the definition immediately, so we omit the proof.

\begin{Pro}\label{OBBPC}
Let $\mathcal{C}$ be an n-extriangulated category. Then a pair $(\mathcal{X},\mathcal{Y})$ of subcategories of $\mathcal{C}$ is a balanced pair if and only if for each object $C\in\mathcal{C}$,
there are the following two distinguished $n$-exangles:
$$\xymatrix{
A\ar[r]&X_{1}\ar[r]&...\ar[r]&X_{n}\ar[r]&C\ar@{.>}[r]^{\delta}&
}$$
$$\xymatrix{
C\ar[r]&Y_{1}\ar[r]&...\ar[r]&Y_{n}\ar[r]&B\ar@{.>}[r]^{\theta}&
}$$
with $X_{n} \in \mathcal{X}$ and $Y_{1}\in \mathcal{Y}$, which are both $\mathcal{C}(\mathcal{X},-)$-exact and $\mathcal{C}(-,\mathcal{Y})$-exact.
\end{Pro}

\begin{Pro}\label{OB-ID}
If $(\mathcal{X},\mathcal{Y})$ is a balanced pair in $\mathcal{C}$, then $(I(\mathcal{X}),I(\mathcal{Y}))$ is an ideal balanced pair.
\end{Pro}

\Pf. It is easy to verify that (BP1) holds. For any object $C\in \mathcal{C}$,
then each induced distinguished $n$-exangle by $\mathcal{X}$-resolution of $C$:
$$\xymatrix{
K_{m}\ar[r]&Z_{m1}\ar[r]&...\ar[r]&Z_{m,n-1}\ar[r]&X_{m}\ar[r]&K_{m-1}\ar@{.>}[r]^{\delta_{m}}&
}$$
such that it is $\mathcal{C}(-,\mathcal{Y})$-exact, where $X_{m}\in \mathcal{X}$.
Thus, for any $y$: $K_{m}\to Y$ with $Y\in \mathcal{Y}$, $y_{\ast}\delta_{m}=0$ by Lemma \ref{lemma1}.
For any $f$: $K_{m}\to N$ in $\mathcal{I}(\mathcal{Y})$, $f$ factors through $M \in \mathcal{Y}$.
i.e., $f=f_{2}f_{1}$, where $f_{1}$: $K_{m}\to M$ and $f_{2}$: $M\to N$.
Consequently, $f_{\ast}\delta_{m}=(f_{2}f_{1})_{\ast}\delta_{m}=(f_{2})_{\ast}(f_{1})_{\ast}\delta_{m}=0$.
Then the condition (BP2) holds. Dually,  the condition (BP3) holds.
\ \hfill $\Box$

\mskip\
Naturally, one may ask whether the converse of Proposition \ref{OB-ID} holds.
To address this question, we introduce the following definition, which is the n-exangulated analogue of Definition 3.2 in \cite{AET}.

\begin{Def}
Let $\mathcal{C}$ be an n-extriangulated category. A negative first extension on $\mathcal{C}$ consists of the following data:

$\mathrm{(1)}$ $\mathbb{E}^{-1}: \mathcal{C}^{\mathrm{op}} \times \mathcal{C} \to \mathrm{Ab}$ is an additive bifunctor.

$\mathrm{(2)}$ For each distinguished $n$-exangle $\xymatrix{A\ar[r]&X_{1}\ar[r]&...\ar[r]&X_{n}\ar[r]&C\ar@{.>}[r]^{\delta}&}$, there exist two natural transformations
$$\xymatrix{
\mathbb{E}^{-1}(-,X_{n-1})\ar[r]&\mathbb{E}^{-1}(-,X_{n})\ar[r]&\mathbb{E}^{-1}(-,C)\ar[r]&\mathcal{C}(-,A)\ar[r]&\cdots\ar[r]&
\mathcal{C}(-,C)
}$$
and
$$\xymatrix{
\mathbb{E}^{-1}(X_{2},-)\ar[r]&\mathbb{E}^{-1}(X_{1},-)\ar[r]&\mathbb{E}^{-1}(A,-)\ar[r]&\mathcal{C}(C,-)\ar[r]&\cdots\ar[r]&\mathcal{C}(A,-).
}$$
\end{Def}

Next, we give the converse of Proposition \ref{OB-ID}.

\begin{Pro}\label{}
If $(\mathcal{I},\mathcal{J})$ is an object ideal balanced pair satisfying $\mathrm{Ob}(\mathcal{I})$ is strongly exact-contravariantly finite and
$\mathrm{Ob}(\mathcal{I})$ is strongly exact-covariantly finite, then $(\mathrm{Ob}(\mathcal{I}),\mathrm{Ob}(\mathcal{J}))$ is a balanced pair.
\end{Pro}

\Pf. By the assumption, for any object $C\in\mathcal{C}$, there exist two distinguished $n$-exangles:
$$\xymatrix{
A\ar[r]^{f_{0}}&X_{1}\ar[r]^{f_{1}}&...\ar[r]&X_{n-1}\ar[r]^{f_{n-1}}&I\ar[r]^{f_{n}}&C\ar@{.>}[r]^{\delta}&
}$$
and
$$\xymatrix{
A\ar[r]^{g_{0}}&J\ar[r]^{g_{1}}&Y_{2}\ar[r]^{g_{2}}&...\ar[r]&Y_{n}\ar[r]^{g_{n}}&B\ar@{.>}[r]^{\delta'}&
}$$
such that $\delta$ is $\mathcal{C}(\mathrm{Ob}(\mathcal{I}),-)$-exact and
$\delta'$ is $\mathcal{C}(-,\mathrm{Ob}(\mathcal{J}))$-exact, where $I\in \mathcal{I}$ and $J\in \mathcal{J}$.
By Lemma \ref{lemma1}, $i^{\ast}\delta=0$ for any $i\in \mathcal{I}$. By Theorem \ref{Th1}, $\delta\in \mathbb{E}_{\mathcal{I}}=\mathbb{E}^{\mathcal{J}}$.
Thus there is a morphism $h_{1}$ such that $g_{0}=h_{1}f_{0}$. By Proposition 3.6 in \cite{HLN}, we have the following commutative diagram:
$$\xymatrix{
A\ar[r]^{f_{0}}\ar@{=}[d]&X_{1}\ar[d]^{h_{1}}\ar[r]^{f_{1}}&...\ar[r]&I\ar[d]\ar[r]^{f_{n}}&C\ar[d]\ar@{.>}[r]^{\delta}&\\
A\ar[r]^{g_{0}}&J\ar[r]^{g_{1}}&...\ar[r]&Y_{n}\ar[r]^{g_{n}}&B\ar@{.>}[r]^{\delta'}&.
}$$
For any $j \in \mathcal{J}$, we consider the following commutative diagram with exact rows:
$$\xymatrix@C=0.5cm{
\mathbb{E}^{-1}(J,M)\ar[r]^{a}\ar[d]&\mathbb{E}^{-1}(A,M)\ar[r]\ar@{=}[d]&\mathcal{C}(B,M)\ar[r]^{b}\ar[d]&\mathcal{C}(Y_{n},M)\ar[r]\ar[d]&
\cdots\ar[r]&\mathcal{C}(J,M)\ar[r]^{c}\ar[d]&\mathcal{C}(A,M)\ar[r]\ar@{=}[d]&\\
\mathbb{E}^{-1}(X_{1},M)\ar[r]^{a'}&\mathbb{E}^{-1}(A,M)\ar[r]&\mathcal{C}(C,M)\ar[r]^{b'}&\mathcal{C}(I,M)\ar[r]&
\cdots\ar[r]&\mathcal{C}(X_{1},M)\ar[r]^{c'}&\mathcal{C}(A,M)\ar[r]&.
}$$
Since $\delta'$ is $\mathcal{C}(-,\mathrm{Ob}(\mathcal{J}))$-exact, $a$ is surjective, and then  $a'$ is also surjective. i.e., $b'$ is injective.
Similarly, we can prove that $c'$ is surjective. i.e., $\delta$ is $\mathcal{C}(-,\mathrm{Ob}(\mathcal{J}))$-exact.
Dually, we obtain that $\delta'$ is $\mathcal{C}(\mathrm{Ob}(\mathcal{I}),-)$-exact.
By Proposition \ref{OBBPC}, $(\mathrm{Ob}(\mathcal{I}),\mathrm{Ob}(\mathcal{J}))$ is a balanced pair.
\ \hfill $\Box$

\subsection{Additive functors}

\begin{Def}
Let $\mathcal{C}$ be an n-exangulated category, and let $\mathbb{E}:\mathcal{C}^{\mathrm{op}} \times \mathcal{C} \to \mathrm{Ab}$ be a biadditive functor.

$\mathrm{(1)}$ A functor $\mathbb{F}: \mathcal{C}^{\mathrm{op}} \times \mathcal{C} \to \mathrm{Set}$ is called a subfunctor of $\mathbb{E}$ if it satisfies the following conditions.

$\bullet$ $\mathbb{F}(C,A)$ is a subset of $\mathbb{E}(C,A)$, for any $A,C\in\mathcal{C}$.

$\bullet$ $\mathbb{F}(c,a) = \mathbb{E}(c,a)\big|_{\mathbb{F}(C,A)}$ holds, for any $a\in\mathcal{C}(A,A')$ and $c\in\mathcal{C}(C',C)$.

In this case, we write $\mathbb{F}\subseteq \mathbb{E}$. $\delta\in\mathbb{F}(C,A)$ is said to be an $\mathbb{F}$-extension.

$\mathrm{(2)}$ A subfunctor $\mathbb{F}\subseteq \mathbb{E}$ is said to be an additive subfunctor if $\mathbb{F}(C,A)\subseteq \mathbb{E}(C,A)$ is an abelian subgroup
for any $A,C\in\mathcal{C}$. In this case, $\mathbb{F}: \mathcal{C}^{\mathrm{op}} \times \mathcal{C} \to \mathrm{Ab}$ itself becomes a biadditive functor.
\end{Def}

Recall that a morphism $i$ is called $\mathbb{F}$-projective if $i^{\ast}\delta=0$ for any $\mathbb{F}$-extension $\delta$.
An additive subfunctor $\mathbb{F}$ is said to have enough projective morphisms if for any $C\in \mathcal{C}$, there exists a distinguished
$n$-exangle $\xymatrix{A\ar[r]&X_{1}\ar[r]&...\ar[r]&X_{n}\ar[r]^{i}&C\ar@{.>}[r]^{\delta}&}$ in $\mathbb{F}$
with $i$ being $\mathbb{F}$-projective.
Dually, we may define $\mathbb{F}$-injective morphisms and the notion that $\mathbb{F}$ has enough injective morphisms.
We denote $\mathbb{F}$-projective morphisms and $\mathbb{F}$-injective morphisms  by $\mathbb{F}$-$\mathbf{proj}$ and $\mathbb{F}$-$\mathbf{inj}$, respectively.
An additive subfunctor $\mathbb{F}$ is called closed under pull-back if $f^{\ast}\delta\in \mathbb{F}$ for any $f$ and $\mathbb{F}$-extension $\delta$.
An additive subfunctor $\mathbb{F}$ is called closed under push-out if $g_{\ast}\theta\in \mathbb{F}$ for any $g$ and $\mathbb{F}$-extension $\theta$.

\begin{Lem}\label{lemma3}
Let $\mathcal{I}$ and $\mathcal{J}$ be two ideals. Then $\mathbb{E}_{\mathcal{I}}$ and $\mathbb{E}^{\mathcal{J}}$ are additive subfunctors of $\mathbb{E}$.
\end{Lem}

\Pf. It is clear that  $\mathbb{E}_{\mathcal{I}}$ contains all split $n$-exangles, and is closed under isomorphisms and direct sums.
We only need to show that $\mathbb{E}_{\mathcal{I}}$ is closed under pull-back and push-out.
Let $\delta\in \mathbb{E}_{\mathcal{I}}(C,A)$. For any $a\in C(A, A')$ and $i\in \mathcal{I}$, $i^{\ast}a_{\ast}\delta=a_{\ast}i^{\ast}\delta=0$. Thus
$a_{\ast}\delta\in \mathbb{E}_{\mathcal{I}}$. That is, $\mathbb{E}_{\mathcal{I}}$ is closed under push-out.
For any $c\in C(C', C)$, $i^{\ast}c^{\ast}\delta=(ci)^{\ast}\delta=0$ since $\mathcal{I}$ is an ideal. Thus
$c^{\ast}\delta\in \mathbb{E}_{\mathcal{I}}$. That is, $\mathbb{E}_{\mathcal{I}}$ is closed under pull-back.
Dually, $\mathbb{E}^{\mathcal{J}}$ is an additive subfunctor.
\ \hfill $\Box$

\begin{Lem}\label{lemma4}
Let $\mathcal{I}$ and $\mathcal{J}$ be two ideals of $\mathcal{C}$. Then the following statements hold.

$\mathrm{(1)}$ $\mathcal{I}$ is strongly precovering if and only if $\mathbb{E}_{\mathcal{I}}$ has enough projective morphisms and $\mathcal{I}=\mathbb{E}_{\mathcal{I}}$-$\mathbf{proj}$.

$\mathrm{(2)}$ $\mathcal{J}$ is strongly preenveloping if and only if $\mathbb{E}^{\mathcal{J}}$ has enough injective morphisms and $\mathcal{J}=\mathbb{E}^{\mathcal{J}}$-$\mathbf{inj}$.
\end{Lem}

\Pf. (1) $(\Rightarrow)$ It is clear that $\mathcal{I}\subseteq\mathbb{E}_{\mathcal{I}}$-$\mathbf{proj}$. For any $C\in \mathcal{C}$, since $\mathcal{I}$ is strongly precovering,
there is a distinguished $n$-exangle:
$$\xymatrix{
A\ar[r]&X_{1}\ar[r] &...\ar[r]&X_{n}\ar[r]^{f_{n}} & C\ar@{.>}[r]^{\delta}&
}$$
with $f_{n}$ being $\mathcal{I}$-precover of $C$. Since $\mathcal{I}\subseteq\mathbb{E}_{\mathcal{I}}$-$\mathbf{proj}$, $f_{n}\in\mathbb{E}_{\mathcal{I}}$-$\mathbf{proj}$.
For any $i$: $X\to C$ in $\mathcal{I}$, $i$ factors through $f_{n}$ since $f_{n}$ is an $\mathcal{I}$-precover. By Lemma \ref{lemma1}, $i^{\ast}\delta=0$. That is, $\delta\in \mathbb{E}_{\mathcal{I}}$.
So $\mathbb{E}_{\mathcal{I}}$ has enough projective morphisms.
For any $g$: $C'\to C$ in $\mathbb{E}_{\mathcal{I}}$-$\mathbf{proj}$, $g$ factors through $f_{n}$ by the definition of $\mathbb{E}_{\mathcal{I}}$-$\mathbf{proj}$.
Since $f_{n}\in \mathcal{I}$, $g\in \mathcal{I}$ as $\mathcal{I}$ is an ideal. That is, $\mathbb{E}_{\mathcal{I}}$-$\mathbf{proj}\subseteq\mathcal{I}$.
Consequently, we have $\mathcal{I}=\mathbb{E}_{\mathcal{I}}$-$\mathbf{proj}$.

$(\Leftarrow)$ For any $M\in \mathcal{C}$, by the assumption, there is a distinguished $n$-exangle:
$$\xymatrix{
N\ar[r]&Y_{1}\ar[r] &...\ar[r]&Y_{n}\ar[r]^{h_{n}} & M\ar@{.>}[r]^{\theta}&
}$$
with $h_{n}\in \mathbb{E}_{\mathcal{I}}$-$\mathbf{proj}=\mathcal{I}$ and $\theta\in \mathbb{E}_{\mathcal{I}}$.
It is sufficient to prove that $h_{n}$ is an $\mathcal{I}$-precover. For any  $i\in\mathcal{I}$, $i^{\ast}\theta=0$ as $\theta\in \mathbb{E}_{\mathcal{I}}$.
By Lemma \ref{lemma1}, $i$ factors through $h_{n}$. That is, $h_{n}$ is an $\mathcal{I}$-precover of $M$.

(2) It is the dual of (1).
\ \hfill $\Box$

\mskip\
The following is the main result of this section.

\begin{Th}
Let $\mathcal{I}$ and $\mathcal{J}$ be two ideals of $\mathcal{C}$. Then the following are equivalent:

$\mathrm{(1)}$ The pair $(\mathcal{I},\mathcal{J})$ is an ideal balanced pair in $\mathcal{C}$.

$\mathrm{(2)}$ There exists an additive subfunctor $\mathbb{F}=\mathbb{E}_{\mathcal{I}}=\mathbb{E}^{\mathcal{J}}$ with enough projective morphisms and enough injective morphisms, and $\mathcal{I}=\mathbb{F}$-$\mathbf{proj}$, $\mathcal{J}=\mathbb{F}$-$\mathbf{inj}$.

\end{Th}

\Pf.  By Theorem \ref{Th1}, Lemma \ref{lemma3} and Lemma \ref{lemma4}, we immediately know that this result holds.
\ \hfill $\Box$

\section{Declarations}

\textbf{Conflict of interest} This work does not have any conflicts of interest.

{\small

}


\begin{thebibliography}{17}

\bibitem{AET} T. Adachi, H. Enomoto and M. Tsukamoto, Intervals of s-torsion pairs in extriangulated categories with negative first extensions,
              Mathematical Proceedings of the Cambridge Philosophical Society, 2023, 174: 451-469.

\bibitem{AR} M. Auslander, I. Reiten, Applications of contravariantly finite subcategories, Advances in Mathematics, 1991, 86(1): 111-152.

\bibitem{AS} J. Asadollahi, S. Sadeghi, Higher ideal approximation theory, Transactions of the American Mathematical Society, 2022, 375(3): 2113-2145.

\bibitem{CHEN} X. Chen, Homotopy equivalences induced by balanced pairs, Journal of Algebra, 2010, 324: 2718-2731.

\bibitem{FGHT} X. Fu, P. Guil-Asensio, I. Herzog, B. Torrecillas, Ideal approximation theory, Advances in Mathematics, 2013, 244: 750-790.

\bibitem{FH} X. Fu, I. Herzog, Powers of the phantom ideal, Proceedings of the London Mathematical Society, 2016, 112(4): 714-752.

\bibitem{FHH} X. Fu, I. Herzog, J. Hu, H. Zhu, Lattice theoretic properties of approximating ideals, Journal of Pure and Applied Algebra, 2022, 226(7):106986.

\bibitem{GKO} C. Geiss, B. Keller, S. Oppermann, $n$-Angulated categories, Journal Fur Die Reine Und Angewandte Mathematik, 2013, 675: 101-120.

\bibitem{GT} R. G$\ddot{o}$bel, J. Trlifaj, Approximations and Endomorphism Algebras of Modules, Walter de Gruyter, Berlin-New York, 2006.

\bibitem{HHZL} J. He, J. He, P. Zhou, Localization of $n$-exangulated categories, 2022, Arxiv: 2205.07644.


\bibitem{HLN} M. Herschend, Y. Liu, H. Nakaoka, $n$-Exangulated categories(I): Definitions and fundamental properties, Journal of Algebra,  2021, 570: 531-586.

\bibitem{HLN2} M. Herschend, Y. Liu, H. Nakaoka, $n$-Exangulated categories(II): Constructions from n-cluster tilting subcategories, Journal of Algebra, 2022, 594: 636-684.

\bibitem{HZhou} J. He, P. Zhou, $n$-Exact categories arising from $n$-exangulated categories, 2021, Arxiv: 2109.12954.

\bibitem{Iyama} O. Iyama. Higher-dimensional Auslander-Reiten theory on maximal orthogonal subcategories, Advances in Mathematics, 2007, 210(1): 22-50.

\bibitem{J} G. Jasso, $n$-Abelian and $n$-exact categories, Mathematische Zeitschrift, 2016, 283: 703-759.

\bibitem{NP} H. Nakaoka and Y. Palu, Extriangulated categories, Hovey twin cotorsion pairs and model structures,
Cahiers de Topologie et Geometrie Differentielle Categoriques, 2019, 60(2): 117-193.

\bibitem{O} F. Ozbek, Precovering and preenveloping ideals, Communications in Algebra, 2015, 43(6): 2568-2584.

\bibitem{WW} Y. Wang, J. Wei, Ideal approximation in $n$-exangulated categories. 2022, arXiv: 2207.13465.

\bibitem{WLH} J. Wang, H. Li, Z. Huang, Applications of exact structures in abelian categories. Publicationes Mathematicae Debrecen, 2016, 88: 269-286.

\bibitem{XF} R. Xu, X. Fu, Ideal balanced pairs in extriangulated categories, Science China Mathematics, 2025, 68(2): 271-284.

\bibitem{XFGS} R. Xu, X. Fu, B.Gao, M. Sun, Ideal approximation theory in extriangulated categories, Journal of Algebra, 2026, 687: 646-676.

\bibitem{ZWei} Q. Zheng, J. Wei. $(n + 2)$-Angulated quotient categories, Algebra Colloquium, 2019, 26(4): 689-720.

\bibitem{ZZ1} P. Zhou and B. Zhu. Cluster-tilting subcategories in extriangulated categories. Theory and Applications of Categories, 2019, 34(8): 221-242.

\bibitem{ZZhou} P. Zhou, B. Zhu, $n$-Abelian quotient categories, Journal of Algebra, 2019, 527: 264-279.





\end{thebibliography}
\end{document}